\documentclass[10pt]{article}
\usepackage{amsfonts}
\usepackage{epsfig,amssymb,amsmath}
\usepackage{graphicx}
\usepackage{subfigure}
\usepackage{multirow}

\newcommand{\calS}{\mathcal{S}}
\newcommand{\calX}{\mathcal{X}}
\newcommand{\calP}{\mathcal{P}}
\newcommand{\calB}{\mathcal{B}}
\newcommand{\calV}{\mathcal{V}}

\newcommand{\calU}{\mathcal{U}}
\newcommand{\PP}{P}

\newcommand{\la}{\langle}
\newcommand{\ra}{\rangle}
\newcommand{\bx}{\boldsymbol{x}}
\newcommand{\bs}{\boldsymbol{s}}
\newcommand{\bu}{\boldsymbol{u}}
\newcommand{\ba}{\boldsymbol{a}}
\newcommand{\bv}{\boldsymbol{v}}
\newcommand{\bg}{\boldsymbol{g}}
\newcommand{\bpi}{\boldsymbol{\Pi}}
\newcommand{\baru}{\bar{\boldsymbol{u}}}
\newcommand{\barx}{\bar{\boldsymbol{x}}}
\newcommand{\barbx}{\bar{\boldsymbol{x}}}
\newcommand{\bxi}{\boldsymbol{\xi}}

\newcommand{\vsigma}{{\varsigma}}
\newcommand{\bvsigma}{\boldsymbol{\varsigma}}
\newcommand{\barvsigma}{\bar{\boldsymbol{\varsigma}}}
\newcommand{\barbvsigma}{\bar{\boldsymbol{\varsigma}}}
\newcommand{\btau}{\boldsymbol{\tau}}
\newcommand{\half}{\frac{1}{2}}
\newcommand{\bff}{\boldsymbol{f}}
\newcommand{\bb}{\boldsymbol{b}}
\newcommand{\real}{\mathbb{R}}

\newcommand{\G}{\mathbf{G}}
\newcommand{\A}{\mathbf{A}}
\newcommand{\bA}{\mathbf{A}}
\newtheorem{theorem}{Theorem}

\newtheorem{assumption}{Assumption}
\newtheorem{proposition}{Proposition}
\newtheorem{algorithm}{Algorithm}
\newtheorem{remark}{Remark}

\newtheorem{example}{Example}

\newcommand{\eb}{\begin{equation}}
\newcommand{\ee}{\end{equation}}

\begin{document}

\title{Canonical Primal-Dual Method for  Solving Non-convex Minimization Problems}
\author{Changzhi Wu$^\dag$, Chaojie Li$^\dag$, and David Yang Gao$^{\dag,\ddag}$\\
{\small $^\dag$  School of Science, Information, Technology \& Engineering,  University of Ballarat,  VIC 3353}\\
{\small $^\ddag$ Research School of Engineering, The Australian National University, Canberra ACT 2601, Australia
 } }
 \maketitle

\begin{abstract}
A new primal-dual  algorithm is presented for solving
 a class of non-convex minimization problems.
 This algorithm is based on  canonical duality theory such that
 the original non-convex minimization problem is first reformulated as
 a convex-concave saddle point optimization problem, which is then
 solved by a quadratically perturbed primal-dual  method. 
 Numerical examples are illustrated.
  Comparing with the existing results, the proposed algorithm
 can achieve better performance.
\end{abstract}

\noindent \textbf{Subject Class:} 49N15, 49M37, 90C26,
90C20\newline \noindent \textbf{Keywords: }  Global optimization, canonical duality theory, quadratic perturbation,
sensor network optimization.

\section{Problems and motivations}
The nonconvex minimization problem to be studied is proposed as the following:
\begin{equation}
(\calP_o): \;\; \min \left\{ \PP\left( \bx\right) = W(\bx) + \half \la \bx, \bA \bx \ra - \la \bx, \bff \ra
  \; | \; \bx \in \calX_a \right\} , \label{ppo}
\end{equation}
 where $\bx= \{ x_i \} \in
\mathbb{R} ^{n}$  is a decision vector, $\A =\left\{
A_{ij}\right\} \in \mathbb{R}^{n\times n}$ is a given real symmetrical matrix,
$\bff =\{ f_i \} \in \mathbb{R}^{n}$ is a given vector,   $\la *
, * \ra $ denotes a bilinear form in $\real^n \times \real^n$;
the feasible space $\calX_a$ is either $\real^n$ or a subset of $\real^n$ with
linear constraints, such that on which, the nonconvex function
 $W:\calX_a\rightarrow
\real$  is well-defined.

Due to the nonconvexity,  Problem ($\calP_o$) may admit many local minima
and local maxima \cite{GaoBook}. It is not an easy task to
identify or numerically compute its global minimizer.
Therefore,
   many numerical methods  have been developed  in  literature,
   including  the extended  Gauss-Newton method
  (see \cite{GN}),  the proximal method (see \cite{PM}), as well as  the popular
  semi-definite programming (SDP) relaxation (see
\cite{kojima-09}).
Generally speaking,     Gauss-Newton type methods
  are local-based such that  only local optimal
solutions can be expected. To find  global
optimal solution often relies on the branch-and-bound \cite{BAB}
as well as the moment matrix based  SDP relaxation  \cite{Ye-opt-08,MomGlo}.
However, these methods are computationally expensive which can be used for solving mainly
  small or medium size problems.

Canonical duality theory has been used successfully for solving a
large class of global  optimization  problems in both continuous
and discrete systems  \cite{GaoBook,gao-cace09}. The main feature
of this theory is that,  depends on the objective function
$W(\bx)$, the nonconvex/nonsmooth/discrete  primal problems can be
transformed into a unified concave maximization problem over a
convex continuous space, which can be solved easily by using
well-developed convex optimization techniques. This potentially
useful theory was developed from
  Gao and Strang's  original work \cite{GaoStrang}
  where the nonconvex function $W(\bx)$ is the so-called {\em stored energy},
  which is required, by the concept (see \cite{marsd-hugh}, page 8),  to be an objective function.
In physics, a real-valued function  is said to be {\em objective}
if only certain fundamental  rules are satisfied (see  \cite{ball} or Chapter 6 in \cite{GaoBook}).
For example, if  $W(\bx)$ is objective, it should be an invariant under certain
coordinate transformations.
Therefore, instead of the decision variables directly, an objective function usually depends on certain
measure (norm) of $\bx$.
 In this
paper, we shall need only the following weak assumptions for the nonconvex function $W(\bx)$ in $(\calP_o)$.
\begin{assumption}\label{assum}
{\em \begin{verse}
    \item[(A1).] There exits a {\em geometrical operator} $\Lambda(\bx): \calX_a \rightarrow \calV_a \subset \real^m$
  and a strictly convex
function $V: \calV_a \subset \real^m \rightarrow \real$ such that
 \eb
 W(\bx) = V(\Lambda(\bx)) \;\; \forall \bx \in \calX_a  . \label{canonical}
 \ee
  \item[(A2).] The geometrical operator $\Lambda(\bx)$ is a vector-valued quadratic mapping in the form of
 \eb
  \Lambda(\bx) =  \left\{ \half \la \bx, \A_1\bx \ra -
   \la \bx,\bb_1\ra,\cdots, \half \la \bx, \A_m\bx \ra - \la
\bx,\bb_m\ra\right\},
\ee
where $\A_i,i=1,\cdots,m$,  are matrices with appropriate
dimensions, and $\bb_i,i=1,\cdots,m, $ are given vectors.
\end{verse}
}
 \end{assumption}

Actually,  Assumption (A1)  is  the so-called {\em canonical transformation}.
Based on this assumption, the proposed  nonconvex problem
 ($\calP_o$) can be reformulated in the following canonical form:
 \eb
 (\calP): \;\; \min \left\{ \PP\left( \bx\right) = V(\Lambda(\bx))  + \half \la \bx, \bA \bx \ra - \la \bx, \bff \ra
  \; : \; \bx \in \calX_a \right\} . \label{pp}
\end{equation}

 The canonical primal problem $(\calP)$  arises naturally from a wide range of applications in engineering and sciences.
For instance, the canonical function $V(\bxi)$ is simply a quadratic function of $\bxi =\Lambda(\bx)$ in the least squares methods for solving  systems of quadratic equations $\Lambda(\bx) = {\bf d} \in \real^m $ (see \cite{ruan-gao-jiao}),
 chaotic dynamical systems  \cite{ruan-gao-ima},   wireless sensor network localization
\cite{gao-runa-pardalos}, general Euclidean distance geometry \cite{Struct},
and computational biology \cite{zhang-gao}.
In computational physics and networks optimization, the   position variable
 $\bx$ is usually a matrix (second-order tensor) and the geometrical operator $\bxi=\Lambda(\bx)$ is
 a positive semi-definite (discredited Cauchy-Riemannian measure)  tensor (see \cite{gao-runa-pardalos}),
 the convex function $V(\bxi)$ is then  an objective function, which is the instance studied by  Gao and Strang  \cite{GaoBook,GaoStrang}.
Particularly,  if  $W(\bx)$ is a quadratic function, the canonical dual problem is equivalent to a SDP problem
(see \cite{gao-runa-pardalos}).
By the facts that the
geometrical operator defined in the Assumption (A2) is a general quadratic mapping,
the nonconvex function $W(\bx)$ studied in this paper is not necessary to be ``objective", which
certainly has extensive applications in complex systems.

The rest of this paper is divided into six sections.
The canonical dual problem is formulated in the next section, where,   some existing difficulties are addressed.
The associated  canonical min-max duality theory is discussed in Section 3.
A proximal point method is proposed in Section 4 to solve this canonical  min-max problem.
Section 5 presents some numerical experiments. Applications to
sensor network optimization are illustrated in Section 6.
The paper is ended  by some some concluding remarks.

\section{Canonical duality theory}

By Assumption (A1),
 the canonical function  $V(\cdot)$ is strictly convex and differentiable on $\calV_a$,
 therefore, the canonical dual mapping $\bvsigma = \nabla V(\bxi) : \calV_a \times \calV_a^* \subset \real^m$ is
 reversible such that  the
following canonical duality relations hold on $\calV_a \times \calV^*_a$.
\begin{equation}
\bvsigma = \nabla V(\bxi) \Leftrightarrow \bxi =
V^{\ast}(\bvsigma) \Leftrightarrow V(\bxi) + V^{\ast}(\bvsigma) =
\la \bxi,\bvsigma \ra, \label{pridul}
\end{equation}
where $V^{\ast}(\bvsigma)$ is the Legendre conjugate of $V(\bxi)$. Clearly,
we have the inverse Legendre conjugate
\begin{equation}
V(\Lambda(\bx)) =\max  \left\{\la \Lambda(\bx), \bvsigma\ra -
V^{\ast}(\bvsigma) \; | \; \bvsigma \in \calV^*_a \right\}.
\label{MaxV}
\end{equation}
By substituting (\ref{MaxV}) into (\ref{pp}), Problem ($\calP$)
can be equivalently  written as
\begin{equation}
 \min_{\bx}\max_{\bvsigma} \left\{
 \Xi(\bx, \bvsigma) \; |
\; (\bx,\bvsigma) \in \calX_a \times \calV^*_a \right\}, \label{EP}
\end{equation}
where  $\Xi: \calX_a \times \calV^*_a \rightarrow \real$ is
 the {\em total complementary function}  defined by
\begin{eqnarray}
\Xi(\bx,\bvsigma)& =& \la \Lambda(\bx),  \bvsigma \ra -
V^{\ast}(\bvsigma) + \half \la \bx , \A \bx \ra - \la \bx , \bff \ra  \notag\\
& = & \half \la \bx,\G(\bvsigma)\bx\ra - V^{\ast}(\bvsigma)
-\la\bx,\btau(\bvsigma)\ra , \label{Xi}
\end{eqnarray}
in which
 \begin{equation}
 \G(\bvsigma) = \A + \sum_{k=1}^m \varsigma_k \A_k, \label{Gkesi}
 \end{equation}
 and
 \begin{equation}
 \btau(\bvsigma) = \bff + \sum_{k=1}^m  \vsigma_k \bff_k.\label{fkeis}
 \end{equation}
 For a given $\bvsigma \in \calV^*_a$, the criticality condition $\nabla_{\bx} \Xi(\bx, \bvsigma) = 0$
 leads to the following {\em canonical equilibrium equation}
 \eb
 \G(\bvsigma) \bx = \btau(\bvsigma).
 \ee
 Let
 \[
\calS_a = \left\{\bvsigma\in \calV^*_a  \;| \; \exists \; \bx\in \calX_a ,
\text{ such that } \G(\bvsigma) \bx = \btau(\bvsigma)  \right\}
\]
be the  dual feasible space, on which,  the canonical dual function is defined by
\begin{equation}
\PP^d(\bvsigma) = \text{sta}\left\{ \Xi(\bx,\bvsigma) \;| \;
\bx\in\calX_a \right\} = - \half \la \G^\dag(\bvsigma) \btau(\bvsigma) , \btau(\bvsigma) \ra
- V^*(\bvsigma) ,  \label{pd}
\end{equation}
 where $\text{sta}\left\{\;\right\}$ stands for finding  value of the
expression in $\{\;\}$ at its stationary  points,
and
$\G^\dag$ represents the generalized  inverse of $\G$.
Particularly, let
\begin{equation}
\calS_a^+ = \left\{\bvsigma\in\calV^*_a\; | \; \G(\bvsigma) \succeq 0
\right\}, \label{Sa}
\end{equation}
where $\G(\bvsigma)\succeq 0$ means that the matrix $\G(\bvsigma)$
is  positive semi-definite.  Clearly,  the total-complementary function
$\Xi(\bx,\bvsigma)$ is convex-concave on  $
\calX_a\times\calS_a^+$, by which,
 the canonical dual problem can be proposed as the following:
\begin{equation}
(\calP^d): \;\;\;\;\;\;\; \max \{
\PP^d(\bvsigma) \; | \; {\bvsigma\in \calS^+_a} \}. \label{ppd}
\end{equation}
The following result is due to the canonical duality theory.
\begin{theorem}[Gao \cite{gao-cace09}]
 Problem $(\calP^d)$ is canonically dual to $(\calP)$ in the sense that if $\barbvsigma$ is a
 critical solution to  $(\calP^d)$, then the vector
\eb
\barbx = \G^\dag (\barbvsigma)\btau(\barbvsigma)
\ee
is a critical point  to $(\calP)$ and
$P(\barbx)  = P^d(\barbvsigma).$

Moreover, if $\barbvsigma \in \calS^+_a$, then $\barbx$ is a global
minimizer of  $(\calP)$  if and only if $\barbvsigma$ is a global maximizer of $(\calP^d)$, i.e.
\eb
P(\barbx) = \min_{\bx \in \calX_a} P(\bx) \;\;  \Leftrightarrow \;\; \max_{\bvsigma \in  \calS^+_a} P^d(\bvsigma) = P^d(\barbvsigma).
\ee
 \end{theorem}

 This theorem shows that if the canonical dual problem $(\calP^d)$ has a critical solution on $\calS^+_a$, then
 the nonconvex primal problem $(\calP)$ is equivalent to a concave maximization dual problem $(\calP^d)$
without duality gap.
      If we further assume that  $\calX_a = \real^n$ and the optimal solution
$\barvsigma$ to   Problem ($\calP^d$) is an interior point of $\calS^+_a$, i.e.,
  $\G(\barvsigma)\succ 0$ , then the  optimal solution
$\barx$ of Problem ($\calP$) can be obtained uniquely  by $\barx
=\G^{-1} (\barvsigma) \btau(\bvsigma)$ (see \cite{GaoWu}).

 However, our experiences show  that for a class of  ``difficult"  global optimization problems,
 the canonical dual problem has no critical solution in $\calS^+_a$ such that $\G(\barvsigma)\succ
 0$. In this paper, we   propose a computational scheme to
solve the case in which the solution is located on the boundary of
$\calS_a^+$, i.e. the dual solution $\barvsigma$ satisfying
$\G(\barvsigma)\in\partial \calS_a^+$.
To continue, we need an additional  mild assumption:
\begin{description}
    \item[A3] There exists an optimal solution $\barbx $ of Problem
    ($\calP$) such that $\G(\barbvsigma )\succeq 0$, where
    $\barbvsigma = \nabla V(\bxi)|_{\bxi =
    \Lambda(\barbx )}$.
\end{description}

 In fact, Assumption
 (A3) is easily satisfied by many real-world problems.
 To see this, let us first examine the following examples.

 \begin{example} {\em  Suppose that $\calX_a$ is a bounded convex
 polytope subset of $\real^n$. Since $\calX_a$ contains only linear
 constraints,  both  $\calV_a$ and $\calS_a$ are also close and bounded. Let $\chi$
 be the smallest eigenvalue of $\sum_{k=1}^m \varsigma_k \A_k$,
  where $\bvsigma=[\varsigma_1,\cdots,\varsigma_m]^T \in \calS_a$.
  Since $\calS_a$ is bounded, $\chi>-\infty$. Let $\bar{\chi}$ be
 the smallest eigenvalue of $\A$. If $\bar{\chi}+\chi \geq 0$, then
 Assumption (A3) is satisfied\footnote{In fact, Problem ($\calP$) is convex  under the condition  $\bar{\chi}+\chi \geq 0$.
 The proof of this result is similar to that of
Proposition 1 given in \cite{PMNon}.}.}

 \end{example}

This example  shows that if the quadratic function $ \half \la \bx , \A \bx \ra  $ is sufficiently
convex,  the non-convexity of $V(\Lambda(\bx))$ becomes
insignificant. Thus, the combination of them is still convex.
However, this is a special  case in nonconvex  systems.
The following example has a wide applications in network optimization.

\begin{example} {\em  Euclidean distance optimization problem:
\begin{equation}
\min \left\{  \sum_{i,j} \left(\|\bx_i - \bx_j \|^2 -d_{i,j}^2\right)^2 +
\sum_{k} \left(\| \bx_k - \ba_k \|^2 -d_k^2 \right)^2 | \; \bx_i \in \real^{ d} \; \forall i = 1, \dots, n \right\} , \label{EDis}
\end{equation}
where $\bx_i $ is the  location vector in Euclidean space $\real^d$,
$d_{i j}$ and $d_k$ are given distance values, the vectors  $\{\ba_k \}$ are
pre-fixed locations. Problem (\ref{EDis}) has many applications,
such as wireless sensor network localization and molecular design, etc.
 For this nonconvex problem, we can choose $\Lambda(\bx)$ to be the
collection of all $\Lambda_{i j}(\bx) = \|\bx_i - \bx_j \|^2
 $ and $\Lambda_k(\bx) = \|\bx_k - \ba_k \|^2 $. In
this case, $V(\bxi)=\sum_{i,j}( \xi_{i j} -d_{i j}^2)^2 + \sum_k (\xi_k-d_k^2)^2$. If
(\ref{EDis}) has the optimal function value of $0$, then
$\xi_{ij} =  d_{i j}^2$ and $\xi_k  = d_k^2$, where
$\bxi=\Lambda(\barbx )$ and $\barbx $ is an optimal solution
of problem (\ref{EDis}). It is easy to check that the dual
variable $\barbvsigma =0 $. Thus, $\det \G(\barbvsigma )=0$.
Therefore, Assumption (A3) holds.
}
 \end{example}

This example shows that
Assumption (A3) is satisfied in the least squares method  for solving
a large class of nonlinear  systems \cite{ruan-gao-ima,ruan-gao-jiao}.
 It is known that  for the conventional SDP
relaxation methods, the   solution of problem
(\ref{EDis}) can be exactly recovered if and only if the SDP solution
of Problem (\ref{EDis}) is a relative interior and the optimal
function value of problem (\ref{EDis}) is $0$ \cite{Tseng}.
 If the problem (\ref{EDis}) has more than one solution,
 the conventional SDP relaxation does not produce any solution.
The goal of this paper is to  overcome  this difficulty by
proposing  a canonical primal-dual iterative scheme.

\section{Saddle-point problem}

Based on  Assumption (A1-A3),  the primal problem ($\calP$) is
relaxed to the following canonical saddle point problem:
\begin{equation}
(\calS\calP): \;\;  \min_{\bx}\max_{\bvsigma} \left\{ \Xi(\bx, \bvsigma) =
 \half \la \bx,\G(\bvsigma)\bx\ra - V^{\ast}(\bvsigma)
-\la\bx,\btau(\bvsigma)\ra \; |  \; (\bx,\bvsigma) \in \calX_a
\times \calS_a^+ \right\}. \label{SP}
\end{equation}

Suppose that $(\barx,\barvsigma)$ is a saddle point of Problem
($\calS\calP$). If $\det(\G(\barvsigma))\neq 0$, we call Problem
($\calS\calP$) is non-degenerate. Otherwise, we call it
degenerate.

\subsection{Problem ($\calS\calP$) is non-degenerate}
\begin{theorem} \label{pro1}
Suppose that Problem ($\calS\calP$) is non-degenerate. Then,
$\barx$ is a unique solution of Problem ($\calP$) if and only if
($\barx,\barvsigma$) is a solution of Problem ($\calS\calP$).
\end{theorem}
Proof. Suppose that ($\barx,\barvsigma$) is the solution of
Problem ($\calS\calP$).  Since Problem ($\calS\calP$) is
non-degenerate, $\G(\barvsigma)\succ 0$, i.e., $\barvsigma \in
\text{int}\calS_a^+$. Thus,
$\nabla_{\bvsigma}\Xi(\barx,\barvsigma) =\barvsigma -
\Lambda(\barx) = 0$. For any $\bx\in \calX_a$, we have
\[
\min_{\bx\in\calX_a}\PP(\bx) =\min_{\bx\in\calX_a}
\max_{\bvsigma\in\calV^*_a} \Xi(\bx,\bvsigma) =
 \min_{\bx\in\calX_a} \max_{\bvsigma\in\calS_a^+}
\Xi(\bx,\bvsigma) = \Xi(\barx,\barvsigma) = \PP(\barx).
\]
Thus, $\barx$ is the optimal solution of Problem ($\calP$).

On the other hand, we suppose that  $\barx$ is the optimal
solution of Problem ($\calP$). Let $\barvsigma = \nabla
V(\Lambda(\barx)) $. Then,
\[
\PP(\barx) = \Xi(\barx,\barvsigma)= \max_{\bvsigma\in\real^m}
\Xi(\barx,\bvsigma).
\]
Since $V(\cdot)$ is strictly convex, we have
\begin{equation}
 \Xi(\barx,\bvsigma) \leq
 \Xi(\barx,\barvsigma) \;\;\forall \;
 \bvsigma\in \calV^*_a \subset \real^m . \label{equ}
\end{equation}
 The equality holds in (\ref{equ}) if and only if $\bvsigma =
\barvsigma$ since $\Xi(\barx,\bvsigma)$ is strictly concave in
terms of $\bvsigma$. Suppose that ($\bx_1,\bvsigma_1$) is also a
saddle point of Problem ($\calS\calP$). By a similar induction as
above, we can show that $\bx_1$ is an optimal solution of Problem
($\calP$). Furthermore, $\PP(\bx_1) = \Xi(\bx_1,\bvsigma_1)$.
Since $\bx_1\in\calX_a$, we have
\[
\PP(\bx_1) = \Xi(\bx_1,\bvsigma_1) \leq \Xi(\barx,\bvsigma_1) \leq
\Xi(\barx,\barvsigma) = \PP(\barx).
\]
The first equality holds only when $\bx_1 = \barx$ since
$\G(\bvsigma_1)\succ 0$. The second inequality becomes equality if
and only if $\bvsigma_1= \barvsigma$ since $V(\cdot)$ is strictly
convex. By the fact that $\barx$ is an optimal solution of Problem
($\calP$) and $\bx_1\in\calX_a$, $\PP(\bx_1) = \PP(\barx)$,
 $\bx_1 = \barx$ and $\bvsigma_1= \barvsigma$. Thus,
 ($\barx,\barvsigma$) is the solution of
Problem ($\calS\calP$). We complete the proof. $\blacksquare$
 
 If $\calX_a = \real^n$, the saddle point
Problem ($\calS\calP$) can be further recast as a convex
semi-definite programming problem.

\begin{proposition}\label{dualsdp}
Suppose that Problem ($\calS\calP$) is non-degenerate and $\calX_a
=\real^n$. Let $\barvsigma$ be the solution of the following
convex SDP problem:
\begin{eqnarray}
(SDP): \;\;\; \min  \left\{  V^{\ast}(\bvsigma) +   g \right\} \;\;\;\;
\text{s.t.}  \;\;
\left[
\begin{array}{cc}
  \G(\bvsigma) & \btau(\bvsigma) \\
  \btau^T(\bvsigma) &  2 g \\
\end{array}
\right] \succeq 0.\label{sdpcost} 
\end{eqnarray}
Then, the SDP
problem defined by (\ref{sdpcost})   has a unique solution ($\bar{g},\barvsigma$) such
that $\G(\barvsigma) \succ 0$. Furthermore, $\barx =
\G^{-1}(\barvsigma)\btau(\barvsigma)$ is the unique solution of
Problem ($\calP$).
\end{proposition}
Proof. By Schur complement lemma \cite{Matrix}, the SDP problem
(\ref{sdpcost})   has a unique solution
($\bar{g},\barvsigma$) such that $\G(\barvsigma) \succ 0$ if and
only if the following convex minimization problem
\begin{eqnarray}
\min \left \{    V^{\ast}(\bvsigma) +
\half \la \G^{-1}(\bvsigma)\btau (\bvsigma), \btau(\bvsigma)\ra | \; \G(\bvsigma) \succeq 0 \right\}  \label{dsdp}
\end{eqnarray}
has a unique solution $\barvsigma$ such that $\G(\barvsigma) \succ
0$. Since $\calX_a = \real^n$, the convex minimization  problem
(\ref{dsdp})  is equivalent to Problem
($\calS\calP$) by  Theorem 3.1 in \cite{GaoWu}.  \hfill $\blacksquare$

\begin{remark}
Theorem \ref{pro1} is actually a special case of the general result obtained by Gao and Strang in finite deformation theory
\cite{GaoStrang}. Indeed,  if we let $\bar{W}(\bx) = W(\bx) + \half \la \bx, \A \bx \ra $ and
 $\bar{\Lambda}(\bx) = \{ \Lambda(\bx), \half \la \bx, \A \bx \ra \}$, then,
 the so-called complementary gap function is simply defined as
 \[
 G (\bx, \bvsigma) = \half \la \bx ,  {\bf G}(\bvsigma) \bx \ra .
 \]
 Clearly, this gap function is strictly positive for any non zero $\bx \in \calX_a$
 if and only if ${\bf G}(\bvsigma) \succ 0 $. Then by Theorem 2 in \cite{GaoStrang} we know that
 the primal problem has a unique solution if the problem ($\calS\calP$) is non-degenerate.
By Theorem \ref{pro1} and Proposition \ref{dualsdp} we know that the nonconvex problem ($\calP$)
 can be solved easily either by
solving a sequence of strict convex-concave
saddle point problems, or via solving a convex semi-definite
programming problem if Problem ($\calS\calP$) is non-degenerate.
By the fact that $g = \half \la \G^{-1}(\bvsigma)\btau (\bvsigma), \btau(\bvsigma)\ra $
is actually the pure complementary gap function (see Eqn (19) in \cite{gao-cace09}),
the convex SDP problem  (\ref{sdpcost}) is indeed
a special case of the canonical dual problem $(\calP^d)$ defined by (\ref{ppd}).
Moreover,   the canonical duality theory can also be used to find the biggest local extrema 
of the nonconvex problem $(\calP)$ (see \cite{GaoWu}).
\end{remark}

\subsection{Problem ($\calS\calP$) is degenerate}
If Problem ($\calS\calP$) is degenerate, i.e.
$\G(\barvsigma)\succeq 0 $ and $\det(\G(\barvsigma)) = 0$ or
$\barvsigma\in\partial \calS_a^+$,  it has multiple saddle points.
The following theorem reveals the relations between Problem
($\calP$) and Problem ($\calS\calP$).

\begin{theorem} \label{th1}
Suppose that Problem ($\calS\calP$) is degenerate.
\begin{description}
    \item[1)] If $\barx$ is a solution of Problem ($\calP$)
    and $\barvsigma = \nabla V(\Lambda(\barx))$, then
    ($\barx,\barvsigma$) is a saddle point of Problem
    ($\calS\calP$).
    \item[2)] If ($\barx,\barvsigma$) is a saddle point of Problem
    ($\calS\calP$), 
    then $\barx$ is a solution of Problem ($\calP$).
     \item[3)] If ($\bx_1,\bvsigma_1$) and ($\bx_2,\bvsigma_2$) are two saddle points of Problem
    ($\calS\calP$), then $\bvsigma_1 = \bvsigma_2$.
\end{description}
\end{theorem}
Proof. 1). Since $\barx$ is a solution of Problem ($\calP$)
    and $\barvsigma = \nabla V(\Lambda(\barx))\in\calS_a^+$ (by Assumption (A3)),
    \[
    \Xi(\barx,\bvsigma) \leq \Xi(\barx,\barvsigma),
    \;\;\forall \bvsigma\in\calS_a^+.
    \]
Furthermore,
\begin{equation}
\langle \nabla \PP(\barx), \bx-\barx \rangle \geq 0,\;\; \forall
\bx\in\calX_a. \label{r1}
\end{equation}
Substituting $\nabla \PP(\barx) = \G(\barvsigma)\barx -
\btau(\barvsigma) = \nabla_{\bx}\Xi(\barx,\barvsigma)$ into
(\ref{r1}), we obtain
\[
\langle \nabla_{\bx}\Xi(\barx,\barvsigma), \bx-\barx \rangle \geq
0,\;\; \forall \; \bx\in\calX_a.
\]
Thus,
\[
\min_{\bx\in\calX_a} \Xi(\bx,\barvsigma) = \Xi(\barx,\barvsigma).
\]
Therefore,
\[
\Xi(\barx,\bvsigma)\leq \Xi(\barx,\barvsigma) \leq
\Xi(\bx,\barvsigma),\;\;\forall\;
(\bx,\bvsigma)\in\calX_a\times\calS_a^+.
\]
This implies that ($\barx,\barvsigma$) is a saddle point of
Problem ($\calS\calP$).

2). Suppose that ($\barx,\barvsigma$) is a saddle point of Problem
($\calS\calP$) and $\nabla_{\bvsigma} \Xi(\barx,\barvsigma) = 0$.
Then,
\[
\PP(\barx) =  \Xi(\barx,\barvsigma) \leq
\Xi(\bx,\barvsigma),\;\;\forall\;
(\bx,\bvsigma)\in\calX_a\times\calS_a^+.
\]
On the other hand,
\[
\Xi(\bx,\barvsigma) = \langle \barvsigma,\Lambda(\bx) \rangle -
V^{\ast}(\barvsigma) - U(\bx) \leq V(\Lambda(\bx)) - U(\bx) =
\PP(\bx).
\]
Combining the above two inequalities,  $\barx$ is a solution
Problem ($\calP$).

3). This result follows directly  from the strict convexity of  both $V(\cdot)$
and $V^{\ast}(\cdot)$. The  proof is completed. \hfill  $\blacksquare$

Theorem \ref{th1} shows that  the nonconvex minimization Problem ($\calP$)
is equivalent to the canonical saddle min-max  Problem ($\calS\calP$).
 What we should emphasize is that the  solutions set of Problem ($\calP$) is in general  nonconvex,
 while   the set
of saddle points of Problem ($\calS\calP$) is  convex.
  For
example, let us consider the following optimization problem:
\begin{equation}
\min \left\{ \half\left((x_1+x_2)^2 -1 \right)^2 +
\half\left((x_1-x_2)^2 -1 \right)^2 | \; \; (x_1, x_2) \in \real^2 \right\} . \label{mulexp}
\end{equation}
Let $\xi=\Lambda(\bx) = [(x_1+x_2)^2 -1, (x_1-x_2)^2 -1]^T $.
Then,
\[
\G(\bvsigma) = \left[
\begin{array}{cc}
  \varsigma_1 + \varsigma_2 & \varsigma_1 - \varsigma_2 \\
  \varsigma_1 - \varsigma_2 & \varsigma_1 + \varsigma_2 \\
\end{array}%
\right],
\]
$V^{\ast}(\bvsigma) = \half \bvsigma^T\bvsigma$. Thus,
$\G(\bvsigma)\succeq 0 \Leftrightarrow \varsigma_1\geq 0 $ and
$\varsigma_2\geq 0$. Clearly, ($\barx,\barvsigma$) is a saddle
point of Problem ($\calS\calP$) if and only if
($\barx,\barvsigma$) is the solution of the following variational
inequality:
\begin{eqnarray}
&&\G(\barvsigma)\barx=0 ,\\
&&\langle \nabla V^{\ast}(\barvsigma) - \Lambda(\barx), \bvsigma -
\barvsigma\rangle \geq 0,\; \forall \bvsigma\geq 0.
\end{eqnarray}
It is easy to verify that the optimization problem (\ref{mulexp})
has four solutions $(1,0)$, $(0,1)$, $(-1,0)$ and $(0,-1)$.
Clearly, its solution set is non-convex. On the other hand, by the statement 3) in
Theorem 3,  we have $\barvsigma = 0$. Thus, $(\barx,\barvsigma)$ is a
saddle point of Problem ($\calS\calP$) if and only if $\barvsigma
= 0$ and $\barx$ satisfies
\begin{eqnarray*}
(x_1+x_2)^2 &\leq& 1,\\
 (x_1-x_2)^2 &\leq& 1.
\end{eqnarray*}
Denote $\Omega = \text{convhull} \{(1,0), (0,1), (-1,0),
(0,-1)\}$, where convhull means convex hull. Therefore, the saddle
point set of Problem ($\calS\calP$) is $\Omega\times 0$ which is a
convex set. This example also shows that the solutions of Problem
($\calP$) are the vertex points of the saddle points set of
Problem ($\calS\calP$).

Now we turn our attention to  the  saddle point  problem
($\calS\calP$).  For some simple optimization
problems, we can simply use   linear perturbation method to solve it. To
illustrate it, let us consider a simple optimization problem given
as below:
\begin{equation}
(\calP_1): \;\;\min_{\bx}\PP_1(\bx) = \half \left(\half\bx^T \A_1
\bx -b_1 \right)^2 + \half \left(\half\bx^T \A_2 \bx -b_2
\right)^2-\la \bx,\bff\ra. \label{p1}
\end{equation}

\begin{proposition} \label{th2}
Suppose that there exists $(\varsigma_1,\varsigma_2)$ such that
$\varsigma_1 \A_1 +\varsigma_2 \A_2 \succ 0$.  If the saddle point
$(\barx,\barvsigma)$ of the associated Problem ($\calS\calP_1$) is
on the boundary of $\calS_a^+$, then for any given $\epsilon>0$,
there exists a $\Delta \bff\in\real^n$ such that $\|\Delta
\bff\|\leq \epsilon$ and the perturbed saddle point Problem
($\calS\calP_1$)
\begin{equation*}
(\calS\calP_1): \;\;  \min_{\bx}\max_{\bvsigma} \left\{ \half \la
\bx,\G(\bvsigma)\bx\ra -\half \bvsigma^T\bvsigma
-\la\bx,\bff+\Delta\bff\ra \; : \; (\bx,\varsigma) \in \real^n
\times \calS_a^+ \right\},
\end{equation*}
has a unique saddle point $(\barx_p,\barvsigma_p)$ such that
$\G(\barvsigma_p)\succ 0$. Furthermore, $\barx_p$ is the unique
solution of
\[
(\calP^{\text{ptb}}_1): \;\;\min_{\bx}\PP_1(\bx) = \half
\left(\half\bx^T \A_1 \bx -b_1 \right)^2 + \half \left(\half\bx^T
\A_2 \bx -b_2 \right)^2 -\la \bx,\bff+\Delta\bff\ra.
\]
where $\G(\bvsigma) = \varsigma_1\A_1 + \varsigma_2\A_2$.
\end{proposition}
Proof. Since $\bx\in\real^n$, Problem ($\calS\calP_1$) is
equivalent to the following optimization problem:
\begin{eqnarray}
\max_{\bvsigma} && - V^{\ast}(\bvsigma) - \half (\bff+\Delta\bff)^T(\bvsigma)\G^{-1}(\bvsigma)(\bff+\Delta\bff) \notag\\
\text{s.t.} && \G(\bvsigma) \succeq 0.\label{expdim2}
\end{eqnarray}
By the assumption that there exists $(\varsigma_1,\varsigma_2)$
such that $\varsigma_1 \A_1 +\varsigma_2 \A_2 \succ 0$, $\A_1$ and
$\A_2$ are simultaneously diagonalizable via congruence. More
specifically, there exists an invertible matrix $C$ such that
\begin{eqnarray*}
\mathbf{C}^T \A_1 \mathbf{C}& =&
\text{diag}(a_1^1,\cdots,a_n^1),\\
\mathbf{C}^T \A_2 \mathbf{C}& =& \text{diag}(a_1^2,\cdots,a_n^2).
\end{eqnarray*}
Under this condition, it is easy to show that for any given
$\epsilon>0$, there exists a $\Delta \bff\in\real^n$ such that
$\|\Delta \bff\|\leq \epsilon$ and
\[
\lim_{\bvsigma\rightarrow \partial \calS_a^+} \half
(\bff+\Delta\bff)^T(\bvsigma)\G^{-1}(\bvsigma)(\bff+\Delta\bff) =
+ \infty.
\]
Thus, the solution of the optimization problem (\ref{expdim2})
cannot be located in the boundary of $\calS_a^+$ for this $\Delta
\bff$. The results follow readily. We complete the proof.
$\blacksquare$

From the  Proposition \ref{dualsdp} we know that  if the solution $\barx$ of
Problem ($\calP_1$) satisfies $\G(\barvsigma)\succ 0$, then it
  can be obtained by simply solving the concave maximization dual problem $(\calP^d)$.
    Otherwise,
Proposition \ref{th2} shows that this solution can be obtained  under a small
perturbation. Thus, the non-convex optimization problem
($\calP_1$) can be completely solved by either the  convex SDP or
the canonical duality. However, for general optimization problems,
the linear perturbation method   may not produce
an interior saddle point of Problem ($\calS\calP$). To overcome
this difficulty, we shall introduce a nonlinear perturbation method  in the next section.

\section{Quadratic Perturbation Method}

We now  focus on the solution of Problem
($\calS\calP$) when it is degenerate. Clearly, Problem
($\calS\calP$) is strictly concave with respect to $\bvsigma$.
However, if Problem ($\calS\calP$) is degenerate, i.e.,
$\barvsigma\in
\partial \calS_a^+$, then Problem ($\calS\calP$) is convex
but not strictly  in terms of $\bx$ . In this  case,
Problem ($\calS\calP$) have multiple solutions. To stabilize such
kind of optimization problems, nonlinear perturbation methods can be used (see \cite{gao-ruan-jogo10}).
Thus, by using the
quadratic perturbation method to Problem ($\calS\calP$),
a  regularized   saddle point problem
 can be proposed as
 \begin{equation}
 \min_{\bx }\max_{\bvsigma\in\calS_a^+}\Xi_{\rho_k}(\bx,\bvsigma)
 = \Xi(\bx,\bvsigma)  + \frac{\rho_k}{2}\|\bx -\bx_k\|^2, \label{regu}
 \end{equation}
where both $\bx_k$ and $\rho_k$, $k=1,2,\cdots,$ are given. In
practical computation, the canonical dual feasible space $\calS_a^+ $
can also be relaxed as
\[
\calS_{\mu_k}^+=\{ \bvsigma\in \calV^*_a \subset \real^m\;:\; \G(\bvsigma)+\mu_k
I\succeq 0\},
\]
 where $\mu_k<\rho_k$. Note that
\[
\Xi_{\rho_k}(\bx,\bvsigma) = \half\la \bx,(\G(\bvsigma)+\rho_k
I)\bx\ra -V^{\ast}(\bvsigma) - \la \bx,\rho_k\bx_k+\btau(\bvsigma)
\ra +\frac{\rho_k}{2}\la \bx_k,\bx_k\ra.
\]
Thus, $\Xi_{\rho_k}(\bx,\bvsigma)$ is strictly convex-concave in
$\real^n\times \calS_{\mu_k}^+$ and
\[
 \min_{\bx
}\max_{\bvsigma\in\calS_{\mu_k}^+}\Xi_{\rho_k}(\bx,\bvsigma) =
\max_{\bvsigma\in\calS_{\mu_k}^+}\min_{\bx
}\Xi_{\rho_k}(\bx,\bvsigma).
\]
 For each given
$\bvsigma\in\calS_{\mu_k}^+$, denote
\[
\bx(\bvsigma) = \arg \min_{\bx }\Xi_{\rho_k}(\bx,\bvsigma).
\]
Then, $\bx(\bvsigma) = (\G(\bvsigma)+\rho_k
I)^{-1}(\rho_k\bx_k+\btau(\bvsigma))$. Substituting this
$\bx(\bvsigma)$ into $\Xi_{\rho_k}(\bx,\bvsigma)$, we obtain
 the perturbed canonical dual function
\[
\PP_{\rho_k}^d(\bvsigma) =
-\half \la (\G(\bvsigma)+\rho_k
I)^{-1}(\rho_k\bx_k+\btau(\bvsigma)),  \rho_k\bx_k+\btau(\bvsigma)\ra - V^{\ast}(\bvsigma)
+\frac{\rho_k}{2}\la \bx_k,\bx_k\ra.
\]
Now   our canonical primal-dual  algorithm can be proposed as follows.

\begin{algorithm} \label{ProxAlg}$\;$
\begin{description}
    \item[Step 1] Initialization $\bx_0$, $\rho_0$, $N$ and the error tolerance $\epsilon$. Set
    $k=0$.
    \item[Step 2] Set $\bvsigma_{k+1} = \arg \max_{\bvsigma\in\calS_{\mu_k}^+}\PP_{\rho_k}^d(\bvsigma)$ and
    $\bx_{k+1} =  (\G(\bvsigma_{k+1})+\rho_{k}
I)^{-1}(\rho_k\bx_k+\btau(\bvsigma_{k+1}))$.
    \item[Step 3] If $\|\bvsigma_{k+1} -\bvsigma_{k} \|\leq \epsilon$, stop.
    Otherwise, set $k=k+1$ and go to Step 2.
\end{description}
\end{algorithm}

\begin{theorem}\label{congencetheorem}
Suppose that
\begin{description}
    \item[1)] $\bar{\rho}\geq\rho_k >0$, $\sigma_k =
    \sum_{i=1}^k\rho_i
    \rightarrow +\infty$, $\rho_k\downarrow 0$, $\mu_k \downarrow 0$ and $0<
    \mu_k<\rho_k$;
    \item[2)] For any given $\bx$,
    $\lim_{\|\bvsigma_k\|\rightarrow\infty}\Xi(\bx,\bvsigma_k)=-\infty$;
     \item[3)] The sequence $\{\bx_k \}$ is a bounded;
\end{description}
Then,  there exists a
$(\barx,\barvsigma)\in\real^n\times\calS_a^+$ such that
$\{\bx_k,\bvsigma_k\}\rightarrow (\barx,\barvsigma)$. Furthermore,
$(\barx,\barvsigma)$ is a saddle point of Problem ($\calS\calP$).
\end{theorem}
\textbf{Proof.} Note that $0< \mu_k<\rho_k$, the perturbed total complementary function 
$\Xi_{\rho_k}(\bx,\bvsigma)$ is strictly convex-concave with
respect to $(\bx,\bvsigma)$ in $\real^n\times\calS_{\mu_k}^+$.
Since $(\bx_{k},\bvsigma_k)$ is generated by Algorithm
\ref{ProxAlg}, we have
\begin{equation}
(\bx_{k},\bvsigma_k) = \arg \min_{\bx
}\max_{\bvsigma\in\calS_{\mu_k}^+}\Xi_{\rho_k}(\bx,\bvsigma) =
\arg \min_{\bx }\max_{\bvsigma\in\calS_{\mu_k}^+}\Xi(\bx,\bvsigma)
+ \frac{\rho_{k-1}}{2}\|\bx -\bx_{k-1}\|^2. \label{rhok}
\end{equation}
That is
\[
\Xi_{\rho_k}(\bx_k,\bvsigma) \leq
\Xi_{\rho_k}(\bx_k,\bvsigma_k)\leq
\Xi_{\rho_k}(\bx,\bvsigma_k),\;\;\forall
(\bx,\bvsigma)\in\real^n\times\calS_{\mu_k}^+.
\]
By the fact that  $\mu_k \downarrow 0$ and $\calS_{\mu_k}^+=\{
\bvsigma\in \calV^*_a \;:\; \G(\bvsigma)+\mu_k I\succeq 0\}$, we have
$\calS_{\mu_k}^+\supseteq \calS_{\mu_{k+1}}^+ $ and
$\bigcap_{k}\calS_{\mu_k}^+ = \calS_a^+$.

 To continue, we suppose that
$(\barx,\barvsigma)$ is a saddle point of Problem ($\calS\calP$),
i.e.,
\[
\Xi(\barx,\bvsigma) \leq \Xi(\barx,\barvsigma)\leq
\Xi(\bx,\barvsigma),\;\;\forall
(\bx,\bvsigma)\in\real^n\times\calS_a^+.
\]
Now we adopt the following steps to prove our results.
\begin{description}
    \item[1)] The  sequence $\{\bx_k\}$ is  convergent, i.e., there
    exists a $\barbx$ such that $\bx_k\rightarrow \barx$.

From  (\ref{rhok}), we have
  \begin{equation}
 \Xi_{\rho_{k-1}}(\bx_k,\bvsigma_k) = \Xi(\bx_k,\bvsigma_k)+ \frac{\rho_{k-1}}{2}\|\bx_k
 -\bx_{k-1}\|^2\leq \Xi_{\rho_{k-1}}(\bx_{k-1},\bvsigma_k) =
 \Xi(\bx_{k-1},\bvsigma_k). \label{24}
   \end{equation}
Clearly,
  \begin{equation}
  \Xi(\bx_{k-1},\bvsigma_k) +  \frac{\rho_{k-2}}{2}\|\bx_{k-1}
 -\bx_{k-2}\|^2  = \Xi_{\rho_{k-2}}(\bx_{k-1},\bvsigma_k). \label{25}
 \end{equation}
  Since $\bvsigma_k\in\calS_{\mu_k}^+\subset\calS_{\mu_{k-1}}^+$
  and $(\bx_{k-1},\bvsigma_{k-1})$ is the saddle point of
  $\Xi_{\rho_{k-1}}(\bx,\bvsigma)$ in
  $\real^n\times\calS_{\mu_{k-1}}^+$, we obtain
  \begin{equation}
 \Xi_{\rho_{k-2}}(\bx_{k-1},\bvsigma_k) \leq
 \Xi_{\rho_{k-2}}(\bx_{k-1},\bvsigma_{k-1})=
 \Xi(\bx_{k-1},\bvsigma_{k-1})+  \frac{\rho_{k-2}}{2}\|\bx_{k-1}
 -\bx_{k-2}\|^2.  \label{26}
  \end{equation}
  Combining (\ref{25}) and (\ref{26}), we obtain
  \[
 \Xi(\bx_{k-1},\bvsigma_k) \leq \Xi(\bx_{k-1},\bvsigma_{k-1}).
  \]
  Thus,
   \begin{equation}
  \Xi(\bx_{k},\bvsigma_k) +  \frac{\rho_{k-1}}{2}\|\bx_{k}
 -\bx_{k-1}\|^2  \leq \Xi(\bx_{k-1},\bvsigma_{k-1}). \label{27}
 \end{equation}
 Repeating the above process, we get
  \begin{equation}
  \Xi(\bx_{k},\bvsigma_k) +  \sum_{i=1}^{k-1}\frac{\rho_{i-1}}{2}\|\bx_{i}
 -\bx_{i-1}\|^2  \leq \Xi(\bx_{1},\bvsigma_{1}). \label{28}
 \end{equation}
 On the other hand,
 \begin{eqnarray}
&&\Xi_{\rho_{k-1}}(\bx_{k},\bvsigma_k) =\Xi(\bx_{k},\bvsigma_k) +
\frac{\rho_{k-1}}{2}\|\bx_k
 -\bx_{k-1}\|^2 \notag\\
 && \geq \Xi_{\rho_{k-1}}(\bx_{k},\barvsigma) = \Xi(\bx_{k},\barvsigma) +
\frac{\rho_{k-1}}{2}\|\bx_k
 -\bx_{k-1}\|^2\notag\\
 &&\geq \Xi(\barx,\barvsigma) +
\frac{\rho_{k-1}}{2}\|\bx_k
 -\bx_{k-1}\|^2. \label{29}
 \end{eqnarray}
 Substituting (\ref{29}) into (\ref{28}) gives rise to
 \[
  \Xi(\barx,\barvsigma) +
\sum_{i=1}^{k-2}\frac{\rho_{i-1}}{2}\|\bx_{i}
 -\bx_{i-1}\|^2  \leq \Xi(\bx_{1},\bvsigma_{1}),\;\;\forall \;k\in\mathbb{N}.
 \]
 Since $\{\bx_k\}$ is a bounded sequence, $\sigma_k
    \rightarrow +\infty$ and $\rho_k\downarrow 0$, the  sequence $\bx_k$ is convergent, i.e., there
    exists a $\barx$ such that $\bx_k\rightarrow \barx$.

     \item[2)] The  sequence $\{\bvsigma_k\}$ is  convergent.
     We first show that $\bvsigma_k$ is a bounded sequence.

In a similar argument to the inequality (\ref{27}), we can show
that
 \[
\Xi(\bx_{k+1},\bvsigma_{k+1}) \geq \Xi(\bx_{k+1},\barvsigma)\geq
\Xi(\barx,\barvsigma).
 \]
 On the other hand,
 \begin{eqnarray*}
&& \Xi_{\rho_k}(\bx_{k+1},\bvsigma_{k+1}) =
\Xi(\bx_{k+1},\bvsigma_{k+1}) +  \frac{\rho_{k}}{2}\|\bx_{k+1}
 -\bx_{k}\|^2 \\
 &&  \leq \Xi_{\rho_k}(\barx,\bvsigma_{k+1}) = \Xi(\barx,\bvsigma_{k+1}) +
 \frac{\rho_{k}}{2}\|\barx
 -\bx_{k}\|^2.
 \end{eqnarray*}
 Summing the above inequalities together yields that
\begin{eqnarray*}
&&\Xi(\barx,\barvsigma) -\frac{\bar{\rho}}{2}\|\barx
 -\bx_{k}\|^2\leq \Xi(\barx,\barvsigma) -\frac{\rho_{k}}{2}\|\barx
 -\bx_{k}\|^2\\
&& \leq  \Xi(\bx_{k+1},\barvsigma) - \frac{\rho_{k}}{2}\|\barx
 -\bx_{k}\|^2 \leq  \Xi(\barx,\bvsigma_{k+1}).
\end{eqnarray*}
By Assumption 2) and $\bx^k\rightarrow\barx$, we know that
$\bvsigma_k$ is a bounded sequence.

Now we suppose that there are two subsequences $\{\bvsigma_k^1\}$
and $\{\bvsigma_k^2\}$ of $\{\bvsigma_k\}$ such that
$\{\bvsigma_k^1\}\rightarrow \bvsigma^1$ and
$\{\bvsigma_k^2\}\rightarrow \bvsigma^2$. Denote $\{\bx_k^1\}$ and
$\{\bx_k^2\}$ are two subsequences of  $\{\bx_k\}$ associated with
$\{\bvsigma_k^1\}$ and $\{\bvsigma_k^2\}$. Clearly,
$\bvsigma^1,\bvsigma^2\in\calS_a^+$. Note that
\begin{eqnarray}
&&\Xi(\bx_{k+1}^1,\bvsigma^2) + \frac{\rho_{k}^1}{2}\|\bx_{k+1}^1
 -\bx_{k}^1\|^2 = \Xi_{\rho_k^1}(\bx_{k+1}^1,\bvsigma^2) \notag\\
 &&\leq
\Xi_{\rho_k^1}(\bx_{k+1}^1,\bvsigma_{k+1}^1)
 = \Xi(\bx_{k+1}^1,\bvsigma_{k+1}^1) + \frac{\rho_{k}^1}{2}\|\bx_{k+1}^1
 -\bx_{k}^1\|^2 . \label{30}
\end{eqnarray}
Thus,
\[
\Xi(\bx_{k+1}^1,\bvsigma^2)\leq \Xi(\bx_{k+1}^1,\bvsigma_{k+1}^1).
\]
Taking limit on both sides of the above inequality yields to
\[
\Xi(\barx,\bvsigma^2) \leq  \Xi(\barx,\bvsigma^1).
\]
In a similar way, we can show that
\[
\Xi(\barx,\bvsigma^1) \leq  \Xi(\barx,\bvsigma^2).
\]
Therefore,
\[
\Xi(\barx,\bvsigma^1) =  \Xi(\barx,\bvsigma^2)
\]
which implies that $\bvsigma^1 = \bvsigma^2$. Hence,
$\{\bvsigma_k\}$ is a convergent sequence.

 \item[3)] We show that if $\{\bx_k,\bvsigma_k\}\rightarrow
(\barx,\barvsigma)$, then $(\barx,\barvsigma)$ is a saddle point
of Problem ($\calS\calP$).

In a similar argument to 2), it is easy to show that for any
$\bvsigma\in\calS_a^+$, we have
\[
\Xi(\barx,\bvsigma) \leq \Xi(\barx,\barvsigma).
\]
So we only need to show that for any $\bx$,
\begin{equation}
\Xi(\barx,\barvsigma) \leq \Xi(\bx,\barvsigma). \label{31}
\end{equation}
Indeed, by the fact that
\[
\Xi_{\rho_k}(\bx_{k+1},\bvsigma_{k+1})\leq
\Xi_{\rho_k}(\bx,\bvsigma_{k+1}),\;\;\forall \bx.
\]
Passing limit to the above inequality yields to the inequality
(\ref{31}). We complete the proof. $\blacksquare$
\end{description}

In Theorem \ref{congencetheorem}, there are three assumptions.
Assumption 1) is on the selection of the parameters and Assumption
2) is always satisfied for strictly convex functions. Assumption
3) is important to ensure the convergence of Algorithm
\ref{ProxAlg}. In fact, from  our numerical experiments,  we found that
 $\bx_k$ might become unbound for certain cases. Therefore, a modified 
algorithm  for solving Problem ($\calP$) is suggested as the  following.

\begin{algorithm} \label{ProxPrimal} \hfill
\begin{description}
    \item[Step 1] Adopt Algorithm \ref{ProxAlg} to solve Problem
    ($\calS\calP$). Denote the obtained solution as
    $(\barx,\barvsigma)$.
    \item[Step 2] If $\|\Lambda(\barx)- \nabla V^*(\barvsigma ) \|\leq
    \epsilon$, output $\barx$ is a global minimizer of Problem ($\calP$), where $\epsilon$ is the tolerance.  Otherwise,
    a gradient-based optimization method is used to refine Problem ($\calP$)
        with initial condition
    $\barx$.
\end{description}
\end{algorithm}
\begin{remark}
Since Problem ($\calS\calP$) is a convex-concave saddle point
problem, many exact and inexact proximal point methods can be
adapted  \cite{RT,PPA1,PPA2}. 
In fact, solving Problem ($\calS\calP$) is an easy task
since it is essentially a convex optimization problem. However, to
obtain a solution of Problem ($\calP$) from the solution set of
Problem ($\calS\calP$) is a difficult task since the
identification of degenerate indices in the nonlinear
complementarity problem is hard \cite{DNCP}. Unlike the classical
proximal point methods,
 our proposed Algorithm \ref{ProxAlg} is based on a sequence of exterior point 
 approximation. In this case, the  gradient operator
 $[\nabla_{\bx}\Xi(\bx,\bvsigma),-\nabla_{\bvsigma}\Xi(\bx,\bvsigma)]$
 in $\real^n\times\calS_a^+$ is not a monotone operator, but $[\nabla_{\bx}\Xi(\bx,\bvsigma)+\mu_k I,-\nabla_{\bvsigma}\Xi(\bx,\bvsigma)]$ is
 monotone in $\real^n\times\calS_{a}^+$. By the fact that $\bigcap_k
 \calS_{\mu_k}^+ =\calS_a^+$, our algorithm generates a convergent
 sequence and its clustering point is a saddle point of Problem
 ($\calS\calP$) under certain conditions. Since
 $[\nabla_{\bx}\Xi(\bx,\bvsigma),-\nabla_{\bvsigma}\Xi(\bx,\bvsigma)]$
 in $\real^n\times\calS_a^+$ is not monotone for each sub-problem, it
 is natural to approximate an optimal solution of Problem
 ($\calP$) under the perturbation of the regularized term $ \frac{1}{2} \rho_k 
 \|\bx-\bx_k\|^2$. This illustrates why our  perturbed (exterior  penalty-type) 
  algorithm usually produces an optimal solution of
 Problem ($\calP$), while the existing proximal point methods based
 on the interior point algorithm do not.
\end{remark}
\begin{remark}
In our proof of Theorem \ref{congencetheorem}, we require that
$\rho_k\rightarrow 0$. For classical proximal point methods, this
condition was not required. In fact, this condition is adopted for
simple proof that of clustering point $(\barx,\barvsigma)$ of the
sequence $\{\bx_k,\bvsigma_k\}$ being a saddle point of Problem
($\calP$). Our simulations show that $\rho_k\rightarrow 0$ can be
relaxed. Indeed, in our test simulations, we found that the
convergence for the case of $\rho_k$ being chosen as a proper
constant parameter is faster than that one of $\rho_k\rightarrow
0$.
\end{remark}

\section{Numerical experiments}

This section  presents some numerical results by 
proposed canonical primal-dual method. In our simulations, the involved SDP is
solved by YALMIP \cite{yalmip} and SeDuMi \cite{SeDuMi}.

\textbf{Example 4.1.} Let us first consider the optimization
problem (\ref{mulexp}). Taking $\rho_k = \frac{1}{k}$ and $\mu_k =
0.1\rho_k$, the initial condition is randomly generated. Table
\ref{exp1} reports the results obtained by our method.

\begin{table}
\caption{Numerical results for optimization problem
(\ref{mulexp})}
\begin{center}
{\scriptsize
\begin{tabular}{ccccc}
  \hline
 Initial condition   & $\barx$ &$\barvsigma$ & $\PP(\barx)=\half\|\barvsigma-\Lambda(\barx)\|^2 $\\
  \hline
$\left(\begin{array}{c}0.81472369\\0.90579194\end{array}\right)$ & $\left(\begin{array}{c}-1.12001364\times10^{-14}\\1.00004756\end{array}\right)$ & $\left(\begin{array}{c}-3.48372378\\-3.48372376\end{array}\right)\times 10^{-9}$&$0.93735607\times 10^{-8}$ \\
$\left(\begin{array}{c}0.60684258\\0.48598247\end{array}\right)$ & $\left(\begin{array}{c}1.00004756\\5.39453096 \times 10^{-14}\end{array}\right)$ & $\left(\begin{array}{c}-3.48358490\\-3.48358548\end{array}\right)\times 10^{-9}$&$0.93735508\times 10^{-8}$ \\
$\left(\begin{array}{c}-0.61543234\\-0.79193703\end{array}\right)$ & $\left(\begin{array}{c}0.56709252\times 10^{-14}\\-1.00004840\end{array}\right)$ & $\left(\begin{array}{c}-3.48379359\\-3.48379378\end{array}\right)\times 10^{-9}$&$0.93735627\times 10^{-8}$ \\
$\left(\begin{array}{c}-0.92181297\\-0.73820724\end{array}\right)$ & $\left(\begin{array}{c}-1.00004756\\0.12834042\times 10^{-13}\end{array}\right)$ & $\left(\begin{array}{c}-3.48370090\\-3.48370051\end{array}\right)\times 10^{-9}$&$0.93735602\times 10^{-8}$ \\
  \hline
\end{tabular}
}
\end{center}
\label{exp1}
\end{table}

From Table \ref{exp1}, we can see that all the four solutions
$(0,1), (1,0), (0,-1)$ and $(-1,0)$ can be detected by our
algorithm with different (randomly generated) initial conditions. The corresponding
$\G(\barvsigma) \approx 0$, as we shown in Proposition \ref{th2},
can also be solved by perturbation method under any given
tolerance. However, the following optimization problem
\begin{equation}
\min_{\bx} \PP(\bx) =\half\sum_{i=1}^m(\bx^T\A_i\bx-d_i)^2,
\label{lea}
\end{equation}
cannot be solved by perturbation method in general, where $\A_i,
i=1,\cdots,m,$ are randomly generated semi-definite matrix and
$d_i,i=1,\cdots,m,$ are chosen such that the optimal function
value of $\PP(\bx)$ is $0$. In fact, $\G(\barvsigma) = 0$ since
the optimal cost function value of the optimization problem
(\ref{lea}) is $0$. Suppose that $m$ is not too small (for example
$m\geq 20$), for any given small perturbation $\Delta \bff$, the
corresponding saddle point problem ($\calS\calP$) has no solution
($\barx, \barvsigma$) such that $\G(\barvsigma)\succ 0$ by our
numerical experiences. Thus, the linear perturbation method cannot be
applied. Now we use our proposed algorithm to solve (\ref{lea})
with different $\rho_k$ and $\mu_k$. In about $80\%$ cases, our
method can capture a solution of Problem ($\calP$). The
corresponding numerical results are reported in Table
\ref{tablea}.

\begin{table}
\caption{Numerical results for optimization problem (\ref{lea})
after $50$ iterations}
\begin{center}
{\scriptsize
\begin{tabular}{ccccc}
  \hline
 $(n,m)$ &$\begin{array}{c} \PP(\barx) \text { with } \rho_k=1/k \\ \text{ and } \mu_k = 0.1\rho_k \end{array}$& $\begin{array}{c} \PP(\barx) \text { with } \rho_k=0.1 \\ \text{ and } \mu_k = 0.1\rho_k
 \end{array}$ \\
  \hline
$(20,25)$ & $4.67244827\times 10^{-6}$ & $4.44146192\times 10^{-8}$ \\
$(30,35)$ & $2.10227829\times 10^{-5}$ & $0.80404292\times 10^{-5}$ \\
$(40,50)$ & $0.00154861$ & $2.34887665\times 10^{-5}$ \\
$(50,60)$ & $0.00951209$ & $0.00032821$ \\
  \hline
\end{tabular}
}
\end{center}
\label{tablea}
\end{table}
During our numerical computation, we observe that for very few
steps (for example, less than $20$ iterations), the numerical
solution by our method is very close to one  solution of Problem
($\calP$). In fact, for all the cases in Table \ref{tablea}, if we
set $\epsilon = 10^{-4}$, then all the obtained results are
satisfied with $\max_i |\bar{x}_i^{\ast}-x^{true}_i| \leq
\epsilon$, $i=1,\cdots, n$, where
$\bx^{true}=[x^{true}_1,\cdots,x^{true}_n]^T$ is one of exact
optimal solutions of Problem ($\calP$). However, it suffers from
slow convergence. Table \ref{tablea} shows it clearly for the last
two cases.  If a gradient-based optimization method is applied,
then the optimal function value is $\PP(\barx)\approx 10^{-8}$ for
all cases in Table \ref{tablea}.

It is obvious that Problem (\ref{lea}) has at least two solutions
because of its symmetry, i.e., if $\barx$ is its solution, so is
$-\barx$. Thus, classical SDP-based relaxation methods in
\cite{Ye-opt-08,Ye-07-MP,kojima-09} cannot produce an exact
solution. However, our method can produce one at the expense of
iterative computation of a sequence of SDPs in most cases.

\section{Applications to Sensor Networks}

In this section, we  apply our proposed method for   sensor
network localization problems.

Consider $N$ sensors and $M$ anchors, both located in the
$d$-dimensional Euclidean space $\mathbb{R}^d$, where $d$ is $2$
or $3$. Let the locations of $M$ anchor points be given as $a_1$,
$a_2$, $\cdots$, $a_M\in \mathbb{R}^d$. The locations of $N$
sensor points $x_1$, $x_2$, $\cdots$, $x_N\in \mathbb{R}^d$ are to
be determined. Let $N_x$ be a subset of $\{(i,j):1\leq i < j \leq
N\}$ in which the distance between the $i$th and the $j$th sensor
point is given as $d_{ij}$ and $N_a$ be a subset of $\{(i,k):1\leq
i \leq N, 1\leq k \leq M\}$ in which the distance between the
$i$th sensor point and the $k$th anchor point is given as
$e_{ik}$. Then, a sensor network localization problem is to find
vector $x_i\in \mathbb{R}^d$ for all $i=1,2,\cdots,N,$ such that
 \begin{eqnarray}
 \|x_i - x_j\|^2 & = & d_{ij}^2,\;\; \forall (i,j) \in N_x, \label{eqn1}\\
 \|x_i - a_k\|^2 & = & e_{ik}^2,\;\; \forall (i,k) \in N_a.
 \label{eqn2}
 \end{eqnarray}
When the given distances $d_{ij}, (i,j)\in N_x,$ and $e_{ik},
(i,k)\in N_a,$ contain noise, the equalities (\ref{eqn1}) and
(\ref{eqn2}) may become infeasible. Thus, instead of solving
(\ref{eqn1}) and (\ref{eqn2}), we formulate it as a non-convex
optimization as given below:
 \begin{equation}
 \min_{x_1,\cdots,x_N}\sum_{(i,j) \in
 N_x}(\|x_i-x_j\|^2-d_{ij}^2)^2 + \sum_{(i,k) \in
 N_a}(\|x_i-a_k\|^2-e_{ik}^2)^2. \label{fourth-opt}
 \end{equation}
Denote $\bx = [x_1^T,\cdots,x_N^T ]^T\in\mathbb{R}^{dN}$. Then,
(\ref{fourth-opt}) can be rewritten as:
\begin{equation}
 \min_{\bx}\left\{\PP(\bx) = \sum_{ij\in N_x}(\bx^T
\A_{ij}\bx - d_{ij}^2)^2 + \sum_{ik\in N_a}(\bx^T
\mathbf{B}_{ii}\bx -2\bff_{ik}^T \bx - (e_{ik}^2
-\bff_{ik}^T\bff_{ik}) )^2\right\}, \label{primalprob}
\end{equation}
where $\A_{ij} = (\mathbf{E}_i - \mathbf{E}_j) (\mathbf{E}_i -
\mathbf{E}_j)^T$, $\mathbf{B}_{ii} = \mathbf{E}_i \mathbf{E}_i^T$,
\[
\mathbf{E}_i =\left(%
\begin{array}{c}
  0_{d\times d} \\
  \cdots \\
  0_{d\times d} \\
  I_{d\times d} \leftarrow i\\
  0_{d\times d} \\
  \cdots \\
  0_{d\times d} \\
\end{array}%
\right) \text{  and  } \bff_{ik}  =\left(%
\begin{array}{c}
  0_{d} \\
  \cdots \\
  0_{d} \\
  a_k \leftarrow i\\
  0_{d} \\
  \cdots \\
  0_{d} \\
\end{array}%
\right).
\]
As in \cite{Ye-07-MP,Kim09}, the root mean square
distance
\begin{equation}
RMSD=\left(\frac{1}{N}\sum_{i=1}^{N}\|\widehat{x}_{i}-x^{*}_{i}\|_{2}^{2}\right)\nonumber
\end{equation}
 is adopted to measure the accuracy of the locations of
the sensor $i$, $i=1,\cdots,N$, where $\widehat{x}_{i}$ and
$x^*_i$ are the estimated position and true positions,
respectively, $i=1,\cdots,N$. The software package SFSDP
\cite{Kim09} is applied for generating test problems and
comparison. During our simulation,
 all of sensors are placed in [0, 1]$\times$[0, 1] randomly and 4
 anchors are fixed at (0.125,0.125), (0.125,0.875),(0.875,0.125),
 and (0.875,0.875), respectively.

For the conventional SDP relaxation methods, the computed sensor
locations match its true locations if and only if the
corresponding sensor network is uniquely localizable
\cite{Ye-opt-08,Ye-07-MP}.  Thus, if the localized sensor network has
multiple solutions, the conventional SDP relaxation methods
\cite{Ye-07-MP, Kim09} fail to produce a good solution of the
optimization problem defined by (\ref{primalprob}). Let us
consider the following network with multiple solutions:
\\

\begin{figure}[htb]
\centering
\includegraphics[width=8.5cm]{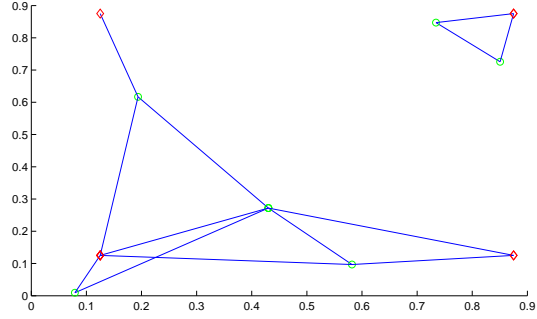}
\caption{Network topology of 6 sensors and 4 anchors}
\label{fig:topo6}
\end{figure}

\textbf{Example 5.1} Consider a sensor network containing 6
sensors and 4 anchors depicted in Figure \ref{fig:topo6}. From
Figure \ref{fig:topo6}, we can see that the sensors $x^{*}_{2}$,
$x^{*}_{3}$ and $x^{*}_{5}$ have two positions.
More specifically, $x_{2}$ can be either (0.0791,0.0091) or
(0.0091,0.1709), $x_{3}$, $x_{5}$ can be either the pair of
$[(0.7342,0.8470), (0.8506,0.7257)]$ or the pair of
$[(1.0158,0.9030), (0.8994,1.0243)]$.  Let $x^{*}$, $\check{x}$
and $\hat{x}$ be the true sensor locations, sensor locations
computed by the SDP method (\cite{kojima-09}), and sensor
locations computed by Algorithm \ref{ProxAlg}, respectively. The
results are depicted in Figure 2 (a) and Figure 2 (c). The true
sensor locations (denoted by circles) and the computed locations
(denoted by stars) are connected by solid lines. From the two
figures, we can clearly see that our method produce better
estimations than the SDP relaxation method in \cite{kojima-09}.
However, we need to solve a sequence of SDPs, but in
\cite{kojima-09}, only one SDP is involved.
\begin{table}
\caption{Numerical results for 6 sensors and 4 anchors}
\begin{center}
{\scriptsize
\begin{tabular}{cccccc}
  \hline
  &True solutions & & Solutions by SDP in \cite{kojima-09}& & Solutions by Algorithm
  \ref{ProxAlg}\\ \hline
$x^{*}_{1}$ & $(0.5818,0.0968)$ & $\check{x}_{1}$ & $(0.5818,0.0961)$ & $\hat{x}_{1}$ & $(0.5818,0.0967)$\\
$x^{*}_{2}$ & $(0.0791,0.0091)$ & $\check{x}_{2}$ & $(0.0775,0.0100)$ & $\hat{x}_{2}$ & $(0.0056,0.1599)$\\
$$ & $(0.0091,0.1709)$ & $$ & $$ & $$\\
$x^{*}_{3}$ & $(0.7342,0.8470)$ & $\check{x}_{3}$ & $(0.7334,0.8985)$ & $\hat{x}_{3}$ & $(1.0167,0.8980)$\\
$$ & $(1.0158,0.9030)$ & $$ & $$ & $$\\
$x^{*}_{4}$ & $(0.1936,0.6169)$ & $\check{x}_{4}$ & $(0.1946,0.6170)$ & $\hat{x}_{4}$ & $(0.1937,0.6169)$\\
$x^{*}_{5}$ & $(0.8506,0.7257)$ & $\check{x}_{5}$ & $(0.7995,0.7439)$ & $\hat{x}_{5}$ & $(0.9047,1.0234)$\\
$$ & $(0.8994,1.0243)$ & $$ & $$ & $$\\
$x^{*}_{6}$ & $(0.4301,0.2720)$ & $\check{x}_{6}$ & $(0.4300,0.2713)$ & $\hat{x}_{6}$ & $(0.4299,0.2718)$\\
  \hline
\end{tabular}
}
\end{center}
\label{tabmul}
\end{table}
To achieve a higher accuracy, we apply the gradient-based
optimization method in SFSDP to refine the solutions obtained by
our method and that obtained by SDP method in \cite{kojima-09}.
After refinement, RMSD obtained by SFSDP is $4.91\times 10^{-5}$
and $2.07\times 10^{-8}$ is obtained by our method. The refined
results are depicted in Figure 2 (b) and (d). From Figure 2 (b),
we observe that there are still big errors for the sensor 3 and
sensor 5 obtained by the refinement of SDP method in
\cite{kojima-09}. Figure 2 (d) shows that our method produces one
of the exact solutions of the optimization problem defined by
(\ref{primalprob}). Thus, our method achieves better performance
no matter before or after refinement.
\begin{figure}[htbp]
\centering \subfigure[Results by SFSDP]{\label{fig:6Na:a}
\begin{minipage}[c]{0.5\textwidth}
\centering
  \includegraphics[height= 6.5cm,width=7cm]{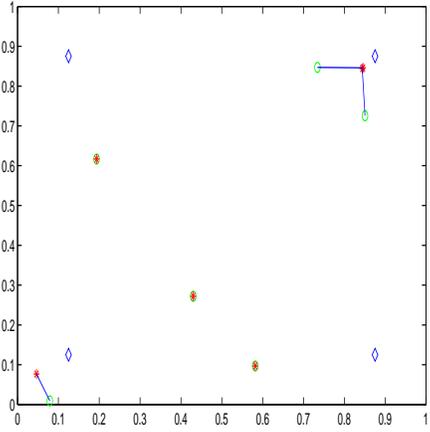}
\end{minipage}%
}%
\subfigure[Results by SFSDP plus the refinement]{
\begin{minipage}[c]{0.5\textwidth}
\centering
  \includegraphics[height= 6.5cm,width=7cm]{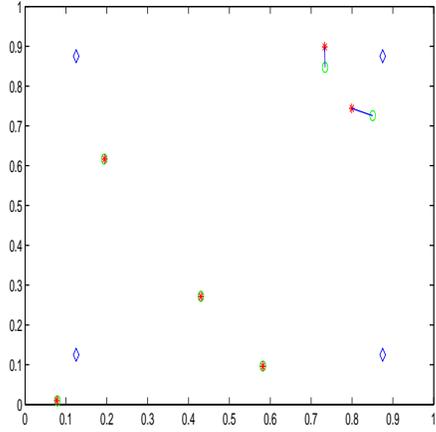}
\end{minipage}
} \subfigure[Results by Algorithm
\ref{ProxAlg}]{\label{fig:6Nb:a}
\begin{minipage}[c]{0.5\textwidth}
\centering
  \includegraphics[height= 6.5cm,width=7cm]{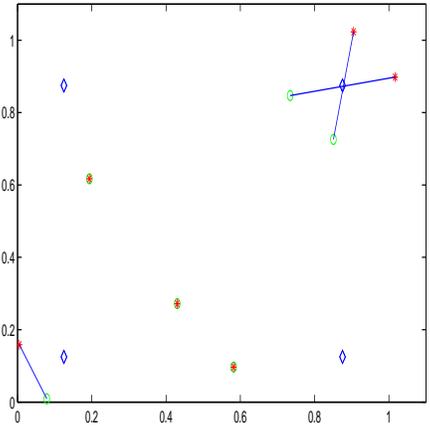}
\end{minipage}%
}%
\subfigure[Results by Algorithm \ref{pridul}]{
\begin{minipage}[c]{0.5\textwidth}
\centering
  \includegraphics[height= 6.5cm,width=7cm]{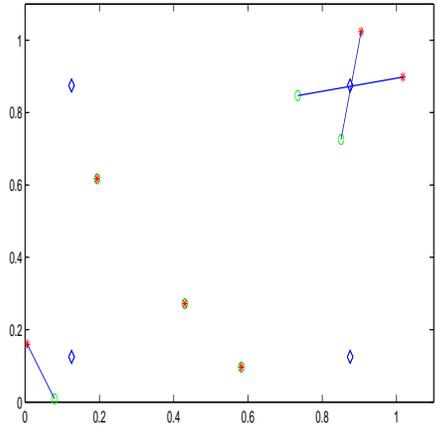}
\end{minipage}
} \caption{Computed locations information of 6 sensors and 4
anchors}\label{fig:6N}
\end{figure}
\\

In practical circumstances, the exact distances $d_{ij}$ and
$e_{ik}$ are unavailable because of the presence of noise during
the measurement. To model such a case, we perturb the distances
as:
\begin{align}
\widehat{d}_{ij}= \max\{(1 + \xi_{ij}), 0.1\}d_{ij} \quad((i, j)\in N_x),\label{34}\\
\widehat{e}_{ik}= \max\{(1 + \xi_{ik}), 0.1\}e_{ik} \quad((i,
k)\in N_a),\label{35}
\end{align}
where $\xi_{ij}$ , $\xi_{ik}$ are random variables and chosen from
the standard normal distribution $N(0, \sigma)$, where $\sigma$ is
the noisy parameter. By substituting (\ref{34}) and (\ref{35})
into (\ref{primalprob}), the corresponding optimization problem
involved in noisy distance is obtained.\\

\textbf{Example 5.2} Consider a sensor network localization
problem with 20 sensors, 4 anchors. Let the radio range be 0.3 and
the noisy parameter be 0.001, respectively. A sensor network
generated randomly by these parameters is depicted in Figure 3.
\begin{figure}[htb]
\centering
\includegraphics[width=8.5cm]{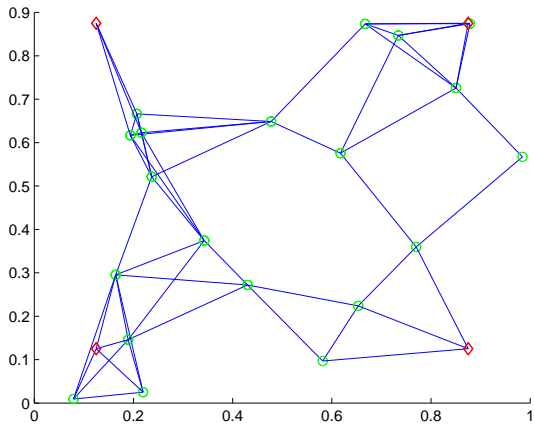}
\caption{Network topology of 20 sensors and 4 anchors}
\label{fig:topo20}
\end{figure}
From Figure 3, we can verify that for this sensor network, it has
a unique solution.

\begin{figure}[htbp]
\centering \subfigure[Results by SFSDP plus the
refinement]{\label{fig:20N:a}
\begin{minipage}[c]{0.5\textwidth}
\centering
  \includegraphics[height= 6.5cm,width=7cm]{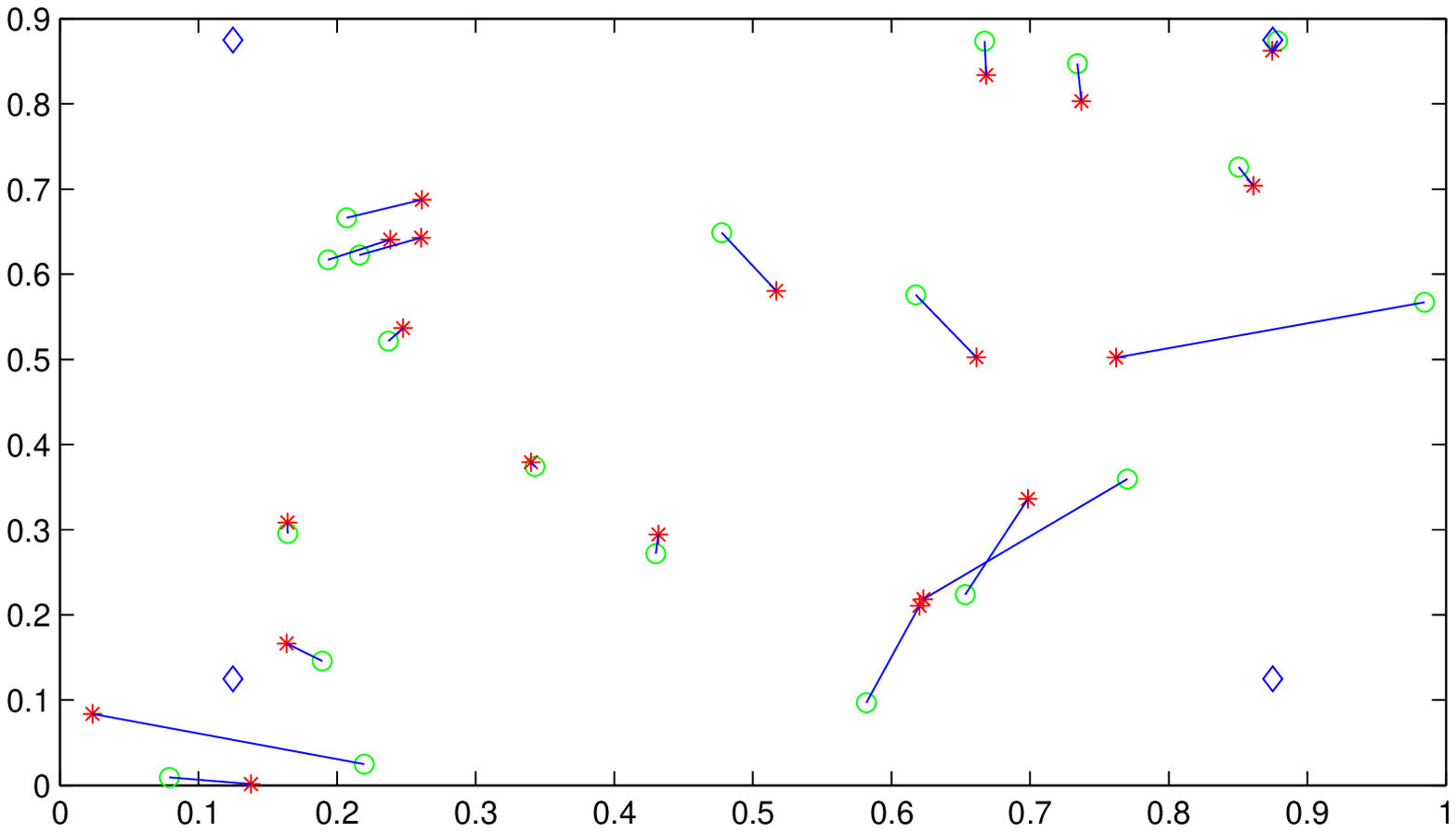}
\end{minipage}%
}%
\subfigure[Results by Algorithm \ref{pridul}]{
\begin{minipage}[c]{0.5\textwidth}
\centering
  \includegraphics[height= 6.5cm,width=7cm]{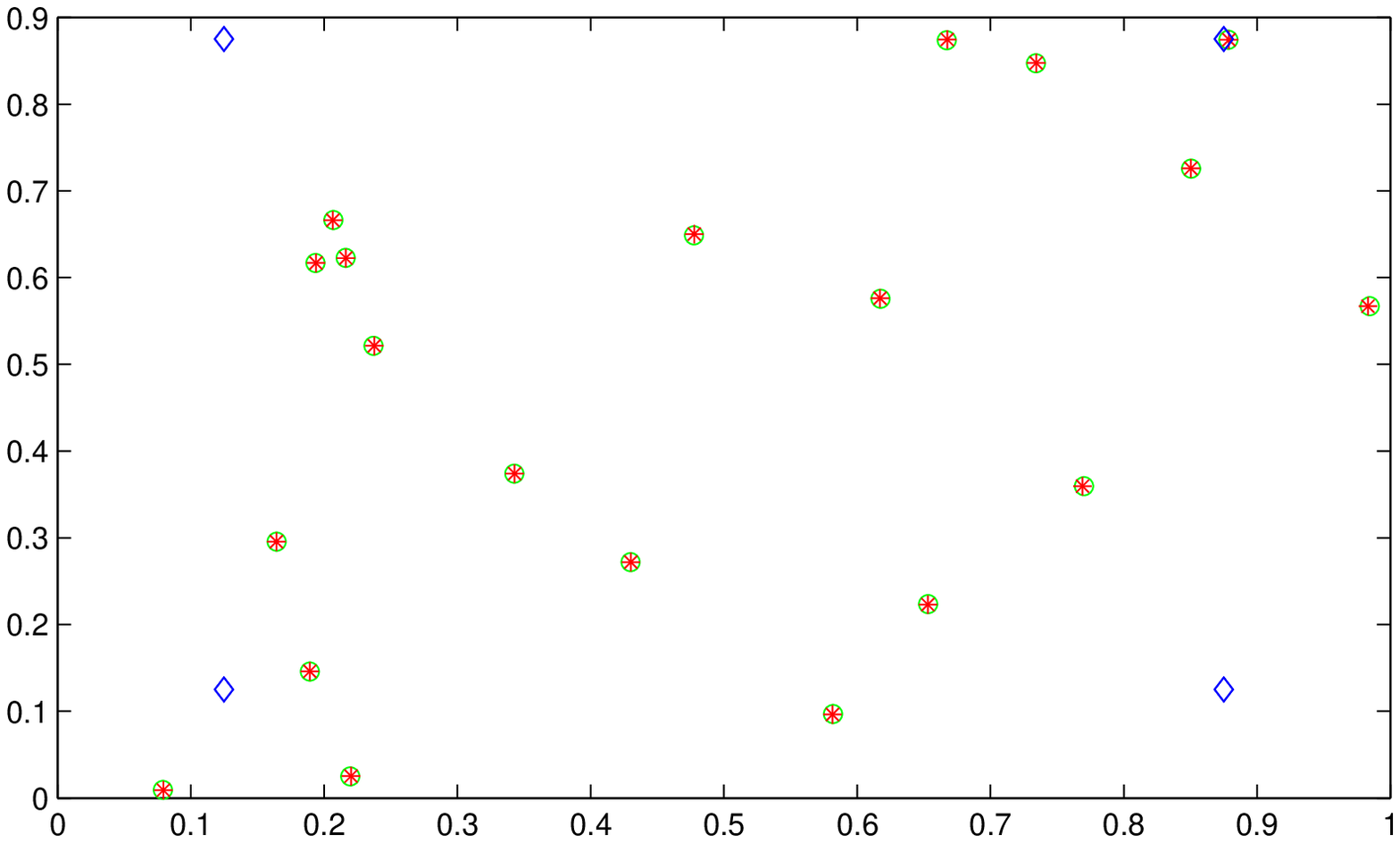}
\end{minipage}
} \caption{Computed locations information of 20 sensors and 4
anchors}\label{fig:20N}
\end{figure}

We apply Algorithm \ref{ProxPrimal} and the SDP method in
\cite{kojima-09} in conjunction with a gradient-based refinement
method to solve it. The computed results are listed in Table
\ref{table20}. The RMSD computed by SFSDP in conjunction with a
gradient-based refinement method is $9.95\times$10$^{-2}$ while
that computed by our method is $4.1041\times 10^{-7}$.  The
computed results by Algorithm \ref{ProxPrimal} and by SDP in
conjunction with a gradient-based refinement method in
\cite{kojima-09} are depicted in Figure \ref{fig:20N}.
From Figure 4 and the values of RMSD, we know that our method
achieves better performance than that by SFSDP in conjunction with
a gradient-based refinement method. This is because if the
distances are inexact, the SDP-based methods in \cite{Ye-07-MP,
kojima-09} are not ensured to produce a good solution. However,
our method is based on the global solution of the optimization
problem defined by (\ref{primalprob}). Thus, the inexact
measurements do not deteriorate the performance of our method.
\begin{table}
\caption{Numerical results for 20 sensors and 4 anchors}
\begin{center}
{\scriptsize
\begin{tabular}{cccccc}
  \hline
  &True solutions & & Solutions by SDP + refinement in \cite{kojima-09}& & Solutions by Algorithm
  \ref{ProxPrimal}\\ \hline
$x^{*}_{1}$ & $(0.5818,0.0968)$ & $\check{x}_{1}$ & $(0.6203,0.2107)$ & $\hat{x}_{1}$ & $(0.5815,0.0963)$\\
$x^{*}_{2}$ & $(0.0791,0.0091)$ & $\check{x}_{2}$ & $(0.1379,0.0015)$ & $\hat{x}_{2}$ & $(0.0795,0.0091)$\\
$x^{*}_{3}$ & $(0.7342,0.8470)$ & $\check{x}_{3}$ & $(0.7369,0.8030)$ & $\hat{x}_{3}$ & $(0.7343,0.8475)$\\
$x^{*}_{4}$ & $(0.1936,0.6169)$ & $\check{x}_{4}$ & $(0.2384,0.6406)$ & $\hat{x}_{4}$ & $(0.1939,0.6168)$\\
$x^{*}_{5}$ & $(0.8506,0.7257)$ & $\check{x}_{5}$ & $(0.8610,0.7040)$ & $\hat{x}_{5}$ & $(0.8503,0.7258)$\\
$x^{*}_{6}$ & $(0.4301,0.2720)$ & $\check{x}_{6}$ & $(0.4319,0.2943)$ & $\hat{x}_{6}$ & $(0.4301,0.2719)$\\
$x^{*}_{7}$ & $(0.9846,0.5671)$ & $\check{x}_{7}$ & $(0.7621,0.5022)$ & $\hat{x}_{7}$ & $(0.9833,0.5670)$\\
$x^{*}_{8}$ & $(0.3429,0.3741)$ & $\check{x}_{8}$ & $(0.3399,0.3793)$ & $\hat{x}_{8}$ & $(0.3430,0.3739)$\\
$x^{*}_{9}$ & $(0.2070,0.6663)$ & $\check{x}_{9}$ & $(0.2612,0.6874)$ & $\hat{x}_{9}$ & $(0.2067,0.6662)$\\
$x^{*}_{10}$ & $(0.6176,0.5756)$ & $\check{x}_{10}$ & $(0.6612,0.5025)$ & $\hat{x}_{10}$ & $(0.6172,0.5762)$\\
$x^{*}_{11}$ & $(0.1644,0.2955)$ & $\check{x}_{11}$ & $(0.1643,0.3085)$ & $\hat{x}_{11}$ & $(0.1643,0.2956)$\\
$x^{*}_{12}$ & $(0.6533,0.2237)$ & $\check{x}_{12}$ & $(0.6984,0.3363)$ & $\hat{x}_{12}$ & $(0.6530,0.2229)$\\
$x^{*}_{13}$ & $(0.6673,0.8736)$ & $\check{x}_{13}$ & $(0.6683,0.8336)$ & $\hat{x}_{13}$ & $(0.6676,0.8746)$\\
$x^{*}_{14}$ & $(0.2161,0.6226)$ & $\check{x}_{14}$ & $(0.2607,0.6429)$ & $\hat{x}_{14}$ & $(0.2165,0.6226)$\\
$x^{*}_{15}$ & $(0.7701,0.3595)$ & $\check{x}_{15}$ & $(0.6232,0.2186)$ & $\hat{x}_{15}$ & $(0.7691,0.3595)$\\
$x^{*}_{16}$ & $(0.1894,0.1458)$ & $\check{x}_{16}$ & $(0.1637,0.1663)$ & $\hat{x}_{16}$ & $(0.1893,0.1460)$\\
$x^{*}_{17}$ & $(0.8786,0.8741)$ & $\check{x}_{17}$ & $(0.8746,0.8626)$ & $\hat{x}_{17}$ & $(0.8789,0.8743)$\\
$x^{*}_{18}$ & $(0.4776,0.6487)$ & $\check{x}_{18}$ & $(0.5169,0.5805)$ & $\hat{x}_{18}$ & $(0.4777,0.6502)$\\
$x^{*}_{19}$ & $(0.2370,0.5215)$ & $\check{x}_{19}$ & $(0.2477,0.5368)$ & $\hat{x}_{19}$ & $(0.2378,0.5215)$\\
$x^{*}_{20}$ & $(0.2197,0.0249)$ & $\check{x}_{20}$ & $(0.0236,0.0836)$ & $\hat{x}_{20}$ & $(0.2202,0.0253)$\\
\hline
\end{tabular}
}
\end{center}
\label{table20}
\end{table}
\\

\begin{figure}[htb]
\centering
\includegraphics[width=8.5cm]{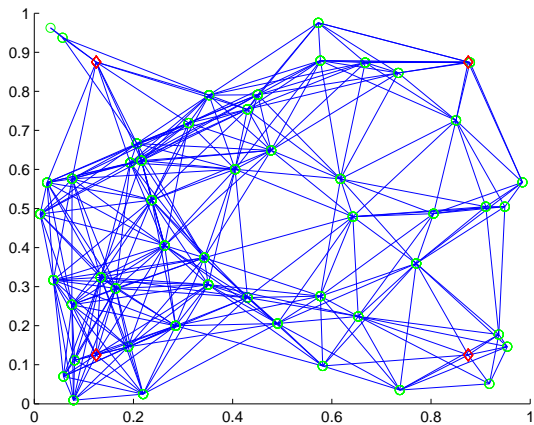}
\caption{Network topology of 50 sensors and 4 anchors}
\label{fig:topo50}
\end{figure}

\textbf{Example 2.3} Consider a sensor network localization
problem with 50 sensors, 4 anchors and noisy perturbation being
0.001. The corresponding connections between sensors and sensors
and sensors and anchors are depicted in Figure \ref{fig:topo50}.

\begin{figure}[htbp]
\centering \subfigure[Results by SFSDP plus the
refinement]{\label{fig:50N:a}
\begin{minipage}[c]{0.5\textwidth}
\centering
  \includegraphics[height= 6.5cm,width=7cm]{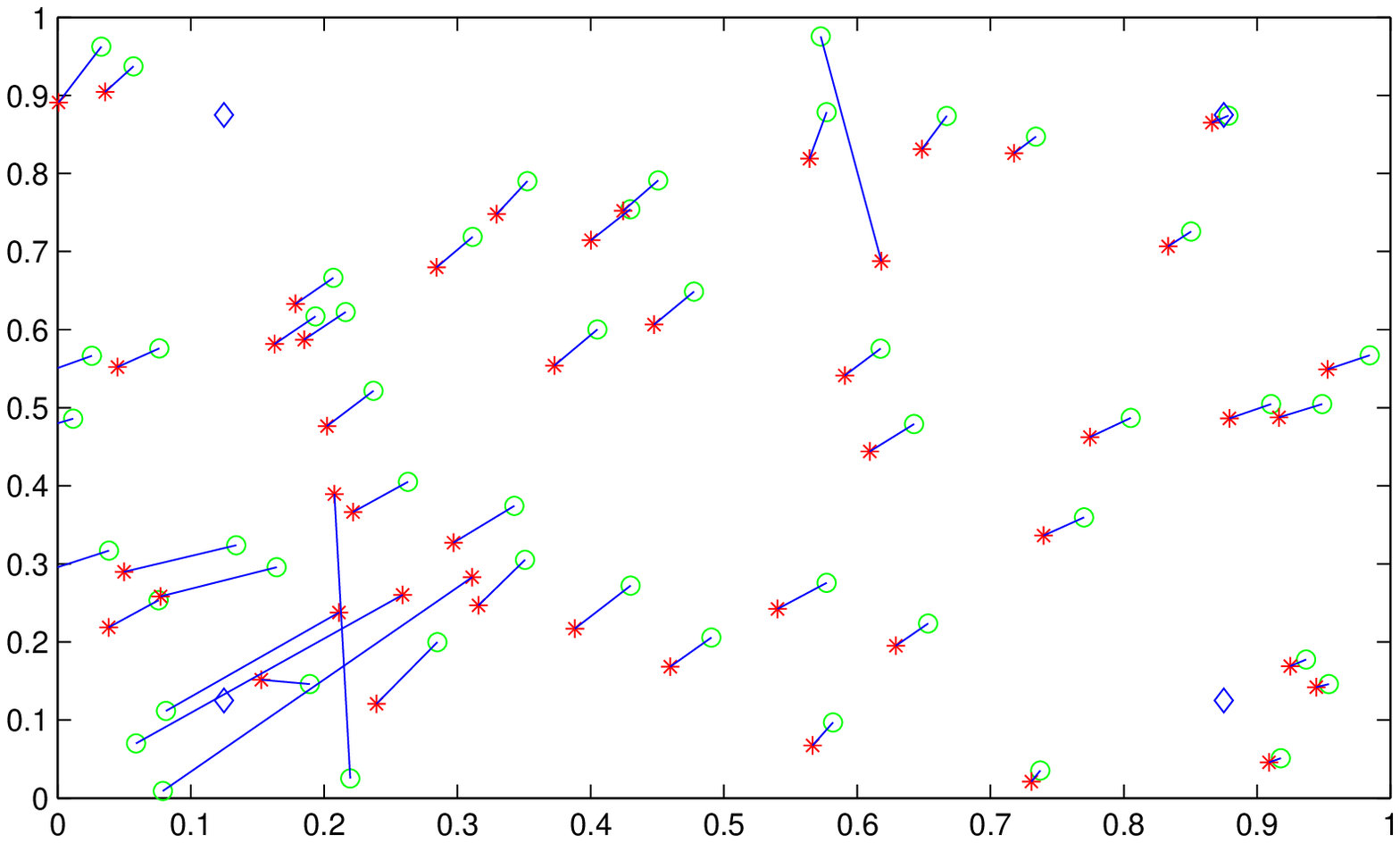}
\end{minipage}%
}%
\subfigure[Results by Algorithm \ref{ProxPrimal}]{\label{fig:50N:b}
\begin{minipage}[c]{0.5\textwidth}
\centering
  \includegraphics[height= 6.5cm,width=7cm]{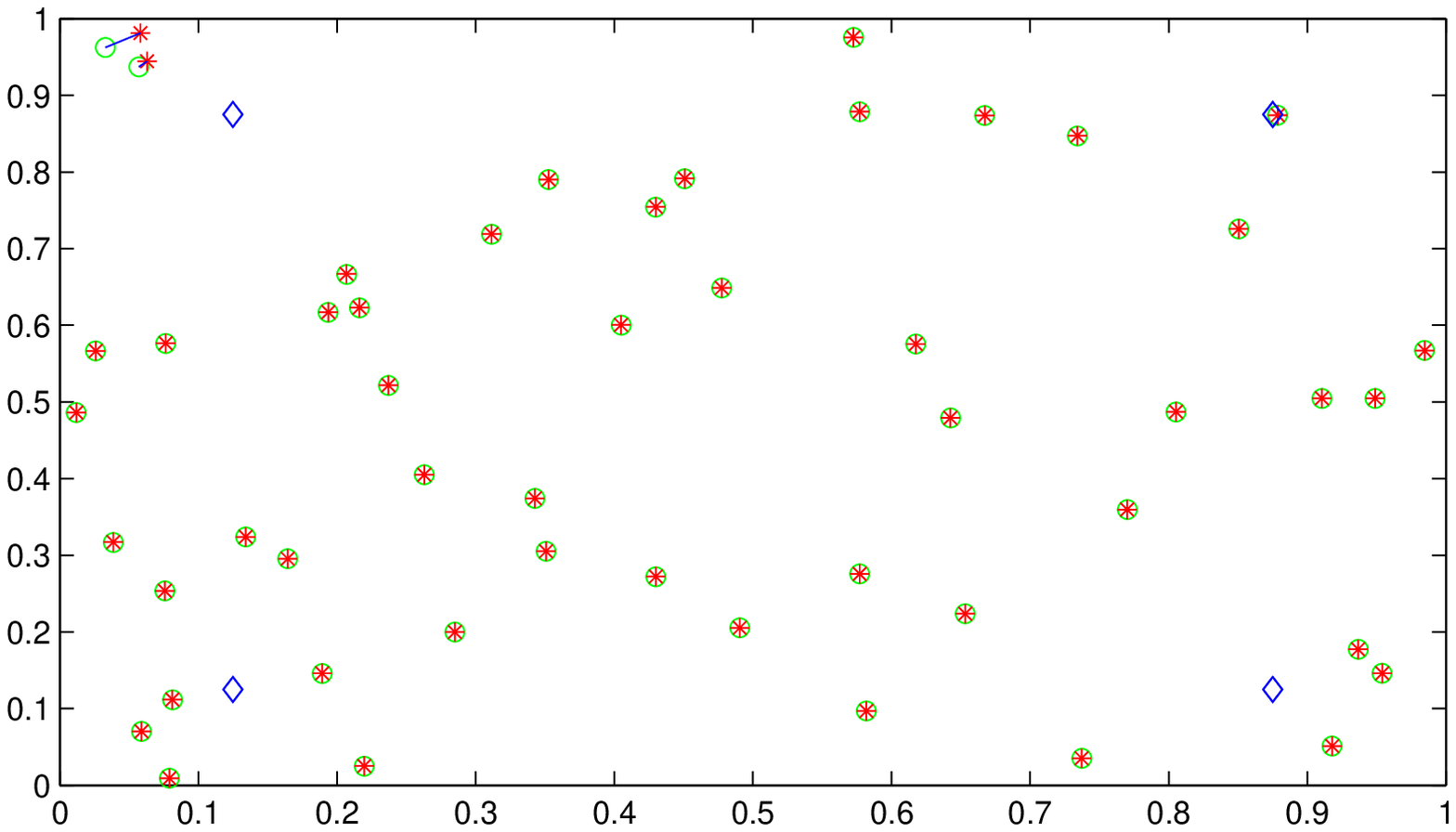}
\end{minipage}
} \caption{Computed locations information of 50 sensors and 4
anchors}\label{fig:50N}
\end{figure}

The computed results by Algorithm \ref{ProxPrimal} and by SFSDP in
conjunction with a gradient-based refinement method is depicted in
Figure \ref{fig:50N}. The RMSD computed by SFSDP in conjunction
with a gradient-based refinement method is $1.07\times 10^{-1}$,
while that by our method is $1.9956\times 10^{-5}$. Both Figure
\ref{fig:50N} and the values of RMSD show that our method achieves
better performance.

\section{Conclusion}
This paper presented an effective method and algorithms for solving a  class of non-convex optimization problems. 
By using the canonical duality theory,
the original non-convex optimization problem is first relaxed to a
convex-concave saddle point optimization problem. Depending on
the singularity of the matrix $\G$,  this  relaxed saddle
point problem is classified in two cases: degenerate or non-degenerate. 
For the non-degenerate case,  the solution of the primal problem can be
recovered exactly through solving a convex SDP problem. Otherwise,  a
quadratic perturbed  primal-dual scheme is proposed  to solve the corresponding degenerate saddle
point problem. We proved that, under certain conditions,  the sequence
generated by our proposed   scheme converges to a solution 
 of the corresponding saddle point problem. If this saddle
point satisfies the condition of 
$\|\Lambda(\barx)- \nabla V^*(\barvsigma ) \|\leq
    \epsilon$ within a given error tolerance,
 then the solution of the primal
problem is also recovered exactly. Otherwise,  $\barx$ is  taken  as a
starting point and   a gradient-based optimization method is applied  to
refine the primal solution. Numerical simulations show that our method can achieve
better performance than the conventional SDP-based relaxation
methods.

\end{document}